\documentclass[final]{siamltex}

\usepackage[dvips]{graphicx}
\usepackage{epsfig}
\usepackage{url}

\usepackage{color}

\xdef\old#1{}
\def\alexo#1{{{ #1}}}

\newtheorem{assumption}[theorem]{\indent {Assumption}}
\newtheorem{algorithm}[theorem]{\indent {Algorithm}}

\def\N{\mathcal{N}}
\def\ox{\overline{x}}

\def\oA{\overline{A}}
\def\A{\mathcal{A}}
\def\oa{\overline{a}}

\def\sin{\sum_{i=1}^n}
\def\sjn{\sum_{j=1}^n}
\def\E{\mathcal{E}}

\begin{document}
\title{Convergence Speed in Distributed Consensus and
Averaging\thanks{This research was supported by the National Science
Foundation under a Graduate Research Fellowship and grants
ECS-0312921 and ECS-0426453. A preliminary version of this paper was presented at the 45th IEEE
Conference on Decision and Control, San Diego, USA, 2006.}}

\author{Alex Olshevsky\thanks{Laboratory for Information and Decision Systems,
Massachusetts Institute of Technology, Cambridge, MA 02139, USA;
{\tt\small alex\_o@mit.edu}} \and John N.\
Tsitsiklis\thanks{Laboratory for Information and Decision Systems,
Massachusetts Institute of Technology, Cambridge, MA 02139, USA;
{\tt\small jnt@mit.edu} }}
 \maketitle

\begin{abstract} We study the convergence speed of distributed
iterative algorithms for the consensus and averaging problems, with
emphasis on the latter. We first consider the case of a fixed
communication topology. We show that a simple adaptation of a
consensus algorithm leads to an averaging algorithm. We prove lower
bounds on the worst-case convergence time for various classes of linear, time-invariant, distributed consensus %averaging
methods, and provide an algorithm that essentially matches those lower bounds.
We then consider the case of a time-varying topology, and
provide a
polynomial-time averaging algorithm. \end{abstract}
\section{Introduction}

Given a set of autonomous agents --- which may be sensors, nodes of
a communication network, cars, or unmanned aerial vehicles
--- the distributed \emph{consensus} problem asks for a distributed
algorithm that the agents can use to agree on an opinion
(represented by a scalar or a vector), starting from different
initial opinions among the agents, and in the presence of possibly severely restricted
communications.

Algorithms that solve the distributed consensus problem provide the
means by which networks of agents can be coordinated. Although each
agent may have access to different local information, the agents can
agree on a decision (e.g., on a common direction of motion, on the
time to execute a move, etc.). Such synchronized behavior %of this sort has
often been observed in biological systems \cite{GVP04}.

The distributed consensus problem has historically appeared in many
diverse areas, such as parallel computation \cite{T84,TBA86,BT89},
control theory \cite{JLM03,OS}, and communication networks
\cite{MSPLM05, LY03}. Recently, the problem has attracted
significant attention \cite{JLM03, LY03, BGPS06, DSW06, MSPLM05,
CSM05, GSMFC05, M05, MVR05, FA05, GCB06, C06, BA06}, motivated by
new contexts and open problems in communications, sensor networks,
and networked control theory. We briefly describe some more of the
more recent applications.

\noindent {\bf Reputation management in ad hoc networks:} It is
often the case that the nodes of a wireless multi-hop network are
not controlled by a single authority or do not have a common
objective. Selfish behavior among nodes (e.g., refusing to forward
traffic meant for others) is possible, and some mechanism is needed
to enforce cooperation. One way to detect selfish behavior is
reputation management: each node forms an opinion by observing the
behavior of its neighbors. One is then faced with the problem of
combining these different opinions into a single globally available
reputation measure for each node. The use of distributed consensus
algorithms for doing this was explored in \cite{LY03}, where a
variation of one of the methods we examine here --- the ``agreement
algorithm'' --- was used as a basis for an empirical investigation.

\noindent {\bf Sensor networks:} A sensor network designed for
detection or estimation needs to combine various measurements into a
decision or into a single estimate. Distributed computation of this
decision/estimate has the advantage of being
  fault-tolerant (network operation is not dependent on a small set of nodes)
and self-organizing (network functionality does not require constant
supervision) \cite{TBA86,BGPS06,DSW06}.

\noindent {\bf Control of autonomous agents:} It is often necessary
to coordinate collections of autonomous agents (e.g., cars or UAVs).
For example, one may wish for the agents to agree on a direction or
speed. Even though the data related to the decision may be
distributed through the network, it is usually desirable that the
final decision depend on all the known data, even though most of them are unavailable at each node. A model motivated by such a context was empirically investigated in \cite{VCBJCS95}.

In this paper, we focus on a special case of the distributed
consensus problem, the distributed {\em averaging} problem.
\alexo{Averaging algorithms guarantee that the final global value
will be the exact average of the initial individual values.} Our
general objective is to characterize the worst-case convergence time
of various averaging algorithms, as a function of the number $n$ of
agents, and to understand their fundamental limitations by providing
lower bounds on the convergence time.

We now outline the remainder of this paper, and preview the main contributions.
In Section \ref{s:2}, we provide some background material, by reviewing the
agreement algorithm of \cite{T84,TBA86} for the distributed
consensus problem. In Sections \ref{s:3}-\ref{s:8}, we consider the case of fixed graphs. In Section \ref{s:3}, we discuss three different ways that the agreement algorithm can provide a solution to the averaging problem.
In particular, we show how an averaging algorithm can be constructed based on two parallel executions of the agreement algorithm.
In Section \ref{s:4},
we define the notions of convergence rate and convergence time, and provide a variational characterization of the convergence rate.

In Section \ref{s:5}, we use results from \cite{LO81} to show that the worst-case convergence time of an averaging algorithm introduced in Section \ref{s:3} is essentially $\Theta(n^3)$.\footnote{Let $f$ and $g$ be two positive functions on the positive integers. We write $f(n)=O(g(n))$ [respectively, $f(n)=\Omega(g(n))$] if there exists a positive constant $c$ and some $n_0$ such that $f(n)\leq c g(n)$ [respectively, $f(n)\geq cg(n)$], for all $n\geq n_0$.
If $f(n)=O(g(n))$ and $f(n)=\Omega(g(n))$ both hold, we write $f(n)=\Theta(g(n))$.
}
  In Section \ref{s:6}, we show that for one of our methods, the convergence rate can be made arbitrarily fast. On the other hand, under an additional restriction that reflects numerical stability considerations, we show that the convergence time of a certain class of algorithms (and by extension of a certain class of averaging algorithms) is $\Omega(n^2)$, in the worst-case. We also provide a simple method (based on executing the agreement algorithm on a spanning tree) whose convergence time essentially matches the
$\Omega(n^2)$ lower bound. In Section \ref{s:7}, we discuss briefly
particular methods that employ doubly stochastic matrices and their
potential drawbacks.

Then, in Section \ref{s:8}, we turn our attention to the case of
dynamic topologies. For the agreement algorithm, we show that its
convergence time for the case of non-symmetric topologies can be
exponentially large in the worst case. On the other hand, for the
case of symmetric topologies, we provide a new averaging algorithm
(and therefore, an agreement algorithm as well), whose convergence
time is $O(n^3)$. To the best our knowledge, none of the existing
consensus or averaging algorithms in the literature, has a similar
guarantee of polynomial time convergence in the presence of
dynamically changing topologies. In Section \ref{s:9}, we report on
some numerical experiments illustrating the advantages of two of our
algorithms. Section \ref{s:10} contains some brief concluding
remarks.

\section{The agreement algorithm\label{s:2}}
The ``agreement algorithm'' is an iterative procedure for the
solution of the distributed consensus problem. It was introduced in
\cite{DeG74} for the time-invariant case, and in \cite{T84,TBA86}
for the case of ``asynchronous'' and time-varying environments. We
briefly review this algorithm and summarize the available
convergence results.

Consider a set $\N=\{1,2,\ldots,n \}$ of nodes. Each node $i$ starts with a scalar value $x_i(0)$;  the vector with
the values of all nodes at time $t$ is %will be
denoted by
$x(t)=(x_1(t),\ldots,x_n(t))$. The agreement algorithm updates
$x(t)$ according to the equation $x(t+1)=A(t)x(t)$, or
$$x_i(t+1)=\sum_{j=1}^n a_{ij}(t)x_j(t),$$ where $A(t)$ is a nonnegative
matrix with entries $a_{ij}(t)$. The row-sums of $A(t)$ are equal to
1, so that $A(t)$ is a stochastic matrix. In particular, $x_i(t+1)$
is a weighted average of the values $x_j(t)$ held by the nodes at
time $t$.

We next state some conditions under which the agreement algorithm is
guaranteed to converge.

\begin{assumption}\label{as:1} There exists a positive constant
$\alpha$ such that:\\ (a) $a_{ii}(t)\geq \alpha$, for all $i$,
$t$.\\ (b) $a_{ij}(t)\in\{0\}\cup [\alpha, 1]$, for all $i$, $j$,
$t$.\\ (c) $\sum_{j=1}^n a_{ij}(t)=1$, for all $i$, $t$.
\end{assumption}

Intuitively, whenever $a_{ij}(t)>0$, node $j$ communicates its
current value $x_j(t)$ to node $i$. Each node $i$ updates its own
value, by forming a weighted average of its own value and the values
it has just received from other nodes. We represent the sequence of
communications between nodes by a sequence
$G(t)=(\N,\mathcal{E}(t))$  of directed graphs, where $(j,i)\in
\mathcal{E}(t)$ if and only if $a_{ij}(t)>0$. Note that $(i,i) \in
\mathcal{E}(t)$ for all $t$, and this condition will remain in
effect throughout the paper.

Our next assumption requires that following an arbitrary time $t$,
and for any $i$, $j$, there is a sequence of communications through
which node $i$ will influence (directly or indirectly) the value
held by node $j$.

\begin{assumption}[Connectivity] \label{as:2} For every $t\geq 0$, the graph
$(\N,\cup_{s\geq t} \mathcal{E}(s))$ is strongly connected.\end{assumption}

Assumption \ref{as:2} by itself is not sufficient to guarantee
consensus (see Exercise 3.1, in p.\ 517 of \cite{BT89}). This
motivates the following stronger version.

\begin{assumption}[Bounded interconnectivity %intercommunication
times]\label{as:3} There is some $B$ such that for all $k$, the
graph $\big(\N,\mathcal{E}(kB)\cup \mathcal{E}(kB+1)\cup\cdots \cup
\mathcal{E}((k+1)B-1)) \big)$ is strongly connected.
\end{assumption}

We note various special cases of possible interest.

\noindent {\em Time-invariant model:} In this model, introduced by
DeGroot \cite{DeG74}, the set of arcs $\mathcal{E}(t)$ is the same
for all $t$; furthermore, the matrix $A(t)$ is the same for all $t$.
In this case, we are dealing with the iteration $x:=Ax$, where $A$
is a stochastic matrix; in particular, $x(t)=A^tx(0)$. Under
Assumptions \ref{as:1} and \ref{as:2}, $A$ is the transition
probability matrix of an irreducible and aperiodic Markov chain.
Thus, $A^t$ converges to a matrix all of whose rows are equal to the
(positive) vector $\pi=(\pi_1,\ldots,\pi_n)$ of steady-state
probabilities of the Markov chain. Accordingly, we have
$\lim_{t\to\infty} x_i(t)=\sum_{i=1}^n \pi_i x_i(0)$.

\noindent {\em Bidirectional model:} In this case, we have $(i,j)\in \mathcal{E}(t)$ if and only if
$(j,i)\in \mathcal{E}(t)$, and we say  that the graph $G$ is \emph{symmetric}.
Intuitively, whenever $i$ communicates to
$j$, there is a simultaneous communication from $j$ to $i$.

\noindent {\em Equal-neighbor model:} Here,
$$a_{ij}(t)=\cases{1/d_i(t),& if $j\in \N_i(t)$,\cr 0,& if $j\notin
\N_i(t)$,\\}$$ where $\N_i(t)=\{j\mid (j,i)\in \mathcal{E}(t)\}$ is
the set of nodes $j$ (including $i$) whose value is taken into
account by $i$ at time $t$, and $d_i(t)$ is its cardinality. This
model is a linear version of a model considered by Vicsek et al.\
\cite{VCBJCS95}. Note that here the constant $\alpha$ of Assumption
\ref{as:1} can be take to be $1/n$.

\begin{theorem}\label{th:1}
Under Assumptions \ref{as:1} and \ref{as:3}, the agreement algorithm
guarantees asymptotic consensus, that is, there exists some $c$
(depending on $x(0)$ and on the sequence of graphs $G(\cdot)$) such
that $\lim_{t\to\infty} x_i(t)=c$, for all $i$.
\end{theorem}

Theorem \ref{th:1} is presented in \cite{TBA86} and proved in
\cite{T84}, in a more general setting that allows for communication
delays, under a slightly stronger version of Assumption \ref{as:3};
see also Ch.~7 of \cite{BT89},  and \cite{TBA86,BHOT05} for
extensions to the cases of communication delays and probabilistic
dropping of packets. The above version of Assumption \ref{as:3} was
introduced in \cite{JLM03}. Under the additional assumption of a
bidirectional model, the bounded interconnectivity time assumption
is unnecessary, as established in \cite{LW04,CMA04} for the
bidirectional equal-neighbor model, and in \cite{HB05, M05} for the
general case.

\section{Averaging with the agreement algorithm in fixed networks}\label{s:3}

In this section, as well as in Sections \ref{s:4}-\ref{s:8}, we assume that the network topology is fixed, i.e., $G(t)=G$ for all $t$, and known. We consider the time-invariant version, $x:=Ax$, of the agreement algorithm, and discuss %introduce
various ways that it can be used to solve the averaging problem. We show that an iteration $x:=Ax$ that solves the consensus problem can be used in a simple manner to provide a solution to the averaging problem as well.

\subsection{Using a doubly stochastic matrix}
As remarked in Section \ref{s:2}, with the time-invariant agreement algorithm $x:=Ax$, we have
\begin{equation}   \lim_{t \rightarrow \infty}
x_i(t)=\sum_{i=1}^n \pi_i x_i(0),\qquad\forall\ i, \label{eq:limit}
\end{equation}
where $\pi_i$ is the steady-state probability \alexo{of node $i$ in}
the Markov chain associated with the stochastic matrix $A$. It
follows that we obtain a solution to the averaging problem if and
only if $\pi_i=1/n$ or every $i$. Since $\pi$ is a left eigenvector
of $A$, with eigenvalue equal to 1, this requirement translates to
the property ${\bf 1}^T A ={\bf 1}^T$, where ${\bf 1}$ is the vector
with all components equal to 1. Equivalently, the matrix $A$ needs
to be doubly stochastic. A particular choice of a doubly stochastic
matrix has been proposed in \cite{OSM04} (see also \cite{GCB06}); it
is discussed further in Sections \ref{s:7} and  \ref{s:9}.

\subsection{The scaled agreement algorithm}\label{sub:scaled}
Suppose that the graph $G$ is fixed a priori and that there is a
system designer or other central authority who chooses a stochastic
matrix $A$ offline, computes the associated steady-state probability
vector (assumed unique and positive), and disseminates the value of
$n\pi_i$ to each node $i$.

Suppose next that the nodes execute the agreement algorithm $\ox:=A\ox$, using the matrix $A$, but with the initial value $x_i(0)$ of each node $i$ replaced by \begin{equation}\label{eq:scale}
\ox_i(0)=
\frac{x_i(0)}{n\pi_i}.
\end{equation}
 Then, the value $\ox_i(t)$ of each node $i$ converges to
$$\lim_{t\to\infty} \ox_i(t)= \sum_{i=1}^n \pi_i \ox_i(0)=
\frac{1}{n}\sum_{i=1}^n x_i(0),$$
and we therefore have a valid averaging algorithm.
This establishes that any (time-invariant) agreement algorithm for the consensus problem translates to an algorithm for the averaging problem as well.
There are two possible drawbacks of the scheme we have just described:
\begin{itemize}
\item[(a)] If some of the $n \pi_i$ are very small, then some of the initial
$\ox_i(0)$ will be very large, which can lead to numerical difficulties \cite{H02}.
\item[(b)]
The algorithm requires some central coordination, in order to choose $A$ and compute $\pi$.
\end{itemize}
The algorithm provided in the next subsection provides a remedy for both of the above drawbacks.

\subsection{Using two parallel passes of the agreement algorithm\label{3a}}

Given a fixed graph $G$, let $A$ be the matrix that corresponds to the
time-invariant, equal-neighbor, bidirectional model
(see Section \ref{s:2} for definitions); in particular, if $(i,j)\in\mathcal{E}$, then $(j,i)\in\mathcal{E}$, %$(i,i)\in\E$,
and
$a_{ij}=1/d_i$, where $d_i$ is the cardinality of $\N_i$.
Under Assumptions \ref{as:1} and \ref{as:2},
the stochastic matrix $A$ is irreducible and aperiodic
(because $a_{ii}>0$ for every $i$). Let $E=\sum_{i=1}^n d_i$.
It is easily verified that the vector $\pi$ with components $\pi_i=d_i/E$, satisfies $\pi^T=\pi^T A$, and is therefore equal to the vector of steady-state probabilities of the associated Markov chain.

The following averaging algorithm employs two parallel runs of the agreement algorithm, with different, but locally determined, initial values.

\begin{algorithm}\label{al:1}
\begin{itemize}
\item[(a)] Each node $i$ sets $y_i(0)=1/d_i$ and $z_i(0)=x_i(0)/d_i$.
\item[(b)] The nodes run the agreement algorithms $y(t+1)=Ay(t)$ and $z(t+1)=Az(t)$.
\item[(c)] Each node sets $x_i(t)=z_i(t)/y_i(t)$.
\end{itemize}
\end{algorithm}

\vspace{5pt}
We have
$$\lim_{t\to\infty} y_i(t)=\sum_{i=1}^n\pi_i y_i(0)=
\sum_{i=1}^n \frac{d_i}{E}\cdot \frac{1}{d_i}=\frac{n}{E},$$
and
$$\lim_{t\to\infty} z_i(t)=\sum_{i=1}^n\pi_i z_i(0)=
\sum_{i=1}^n \frac{d_i}{E}\cdot \frac{x_i(0)}{d_i}=\frac{1}{E}\sum_{i=1}^n x_i(0).$$
This implies that
$$\lim_{t\to\infty} x_i(t)= \frac{1}{n}\sum_{i=1}^n x_i(0),$$
i.e., we have a valid averaging algorithm.
Note that the iteration $y:=Ay$ need not be repeated
if the network remains unchanged and the averaging algorithm is to be executed again with different initial opinions. Finally, if $n$ and $E$ are known by all nodes, the iteration $y:=Ay$ is unnecessary, and we could just set $y_i(t)=n/E$.

\section{Definition of the convergence rate and the convergence time}\label{s:4}

The convergence rate of any of the algorithms discussed in Section \ref{s:3} is determined by the convergence rate of the matrix powers $A^t$. In this section, we give a definition of the convergence rate (and convergence time) and provide a tool for bounding the convergence rate. As should be apparent from the discussion in Section \ref{s:3}, there is no reason to restrict to doubly stochastic matrices, or even to nonnegative matrices. We therefore start by
specifying the class of matrices that we will be interested in.

Consider a matrix $A$ with the following property: for every $x(0)$, the sequence generated by letting $x(t+1)=Ax(t)$ converges to $c{\bf 1}$, for some scalar $c$. Such a matrix corresponds to a legitimate agreement algorithm, and can be employed in the scheme of Section 3.2 to obtain an averaging algorithm, as long as 1 is an eigenvalue of $A$ with multiplicity 1, and the corresponding
left eigenvector, denoted by $\pi$, has nonzero entries.
Because of the above assumed convergence property, all other eigenvalues must have magnitude less than 1. Note, however, that we allow $A$ to have some negative entries.

Suppose that $A$ has the above properties. Let
$1=\lambda_1,\lambda_2,\ldots,\lambda_n$, be the eigenvalues of $A$,
sorted in order of decreasing maginitude. We also let $X$ be the set
of vectors of the form $c{\bf 1}$, i.e., with equal components.
Given such a matrix $A$, we define its {\em convergence rate},
$\rho$, by \begin{equation}
 \rho = \sup_{x(0)
\notin X} \lim_{t \rightarrow \infty} \Big(\frac{\|x(t) -
x^{*}\|_2}{\|x(0)-x^{*}\|_2}\Big)^{1/t}, \label{rhodef}
\end{equation} where  $x^*$ stands for $\lim_{t \rightarrow \infty}
x(t)$.
As is well known, we have $\rho=\max\{|\lambda_2|,|\lambda_n|\}$.

We also define the {\em convergence time} $T_n(\epsilon)$, by
$$T_n(\epsilon) = \min\Big\{
\tau: \frac{\|x(t) - x^*\|_{\infty}}{\|x(0) - x^*\|_{\infty}}
\leq \epsilon,\ \ \forall\ t \geq \tau,\ \forall\ x(0)\notin X\Big\}.$$
Although we use the
infinity norm to define the convergence time, bounds for other norms
can be easily obtained from our subsequent results, using the
equivalence of norms.

\alexo{Under the above assumptions, a result from \cite{XB04} states
\[ \rho = \max \{ |\lambda_2|, |\lambda_n|\}.\] To study $\rho$,
therefore, we must develop techniques to bound the eigenvalues of
the matrix $A$. To this end, } we will be using the following result
from \cite{LO81}. We present here a slightly more general version,
and include a proof, for completeness.

\begin{theorem}\label{th:2} Consider an $n\times n$ matrix $A$ and let
$\lambda_1,\lambda_2,\ldots,\lambda_n$, be its eigenvalues, sorted in order of decreasing maginitude.
Suppose that the following conditions hold.
\begin{itemize}
\item[(a)] We have $\lambda_1=1$ and $A {\bf 1} ={\bf 1}$.
\item[(b)] There exists a positive vector $\pi$ such that $\pi^T A=\pi^T$.
\item[(c)] For every $i$ and $j$, we have $\pi_i a_{ij}=\pi_j a_{ji}$.
\end{itemize}
Let
$$S=\Big\{x\ \Big|\  \sum_{i=1}^n \pi_i x_i =0,\ \sum_{i=1}^n \pi_i x_i^2=1\Big\}$$
Then, all eigenvalues of $A$ are real, and
\begin{equation}\label{eq:gg}
\lambda_2=1-\frac{1}{2} \min_{x\in S}\sum_{i=1}^n \sum_{j=1}^n \pi_i a_{ij}
(x_i-x_j)^2.\end{equation}
In particular, for any vector $y$ that satisfies $\sin \pi_i y_i=0$, we have
\begin{equation} \lambda_2 \geq 1 - \frac{\displaystyle{\sin\sjn \pi_i a_{ij}
 (y_i - y_j)^2}}
{\displaystyle{2\sin \pi_i y_i^2}}.
 \label{mineigensym}
\end{equation}
\end{theorem}
\begin{proof}
Let $D$ be a diagonal matrix whose $i$th diagonal entry
is $\pi_i$. Condition (c) yields $DA=A^TD$. We define the inner product
$\langle\cdot,\cdot\rangle_{\pi}$ by
$\langle x,y\rangle_{\pi}=x^TDy$. We then have
$$\langle x,Ay\rangle_{\pi}=x^TDAy=x^TA^T Dy=\langle Ax,y\rangle_{\pi}.$$
Therefore, $A$ is self-adjoint with respect to this inner product,
which proves that $A$ has real eigenvalues.

Since the
largest eigenvalue is $1$, with an eigenvector of ${\bf 1}$, we use the variational
characterization of the eigenvalues of a self-adjoint matrix
(Chapter 7, Theorem 4.3 of \cite{T05}) to obtain
\begin{eqnarray*}
\lambda_2 &=& \max_{x\in S} \langle x,Ax\rangle_{\pi}\\
 & =
&  \max_{x\in S} \sum_{i=1}^n\sum_{j=1}^n \pi_i a_{ij}
x_i x_j \\
 &=& \frac{1}{2} \max_{x\in S} \sum_{i=1}\sum_{j=1}
\pi_i a_{ij}(x_i^2 + x_j^2 - (x_i - x_j)^2).
\end{eqnarray*}
For $x \in S$, we have
$$\sum_{i=1}^n \sum_{j=1}^n \pi_i a_{ij} (x_i^2 +
x_j^2) =2\sin\sjn \pi_i a_{ij} x_i^2 = 2\sum_{i=1}^n \pi_i x_i^2 =
\alexo{2 \langle x,x \rangle_{\pi}=} 2,$$ which yields
$$
\lambda_2 = 1- \frac{1}{2}\min_{x\in S}\sum_{i=1}^n\sum_{j=1}^n
\pi_i a_{ij} (x_i - x_j)^2. \label{mineigeq} $$ Finally,
Eq.~(\ref{mineigensym}) follows from (\ref{eq:gg}) by considering
the vector $x_i=y_i/\sqrt{(\sjn \pi_j y_j^2)}$.
\end{proof}

Note that the bound of Eq.\ (\ref{mineigensym}) does not change if
we replace the vector $y$ with $\alpha y$, for any $\alpha \neq 0$.

\section{Convergence time for the algorithm of Section \ref{3a}}\label{s:5}

For the equal neighbor, time-invariant, bidirectional model, tight bounds on the
convergence rate were derived in \cite{LO81}.

\begin{theorem}\label{th:3} \cite{LO81}
Consider the equal-neighbor, time-invariant, bidirectional model, on a connected graph with $n$ nodes.
The convergence rate satisfies
\[ \rho \leq 1 - \gamma_1 n^{-3}, \]
where $\gamma_1$ is a constant independent of $n$. Moreover, there
exists some $\gamma_2>0$ such that for every positive integer $n$,
there exists an $n$-node connected symmetric graph for which
\[ \rho \geq 1 - \gamma_2 n^{-3}.
\]
\end{theorem}

Theorem \ref{th:3} is proved in \cite{LO81} for the case of
symmetric graphs without self-arcs. It is not hard to check that
essentially the same proof goes through when self-arcs are present,
the only difference being in the values of the constants $\gamma_1$
and $\gamma_2$. This is intuitive because the effect of the
self-arcs is essentially a ``slowing down'' of the Markov chain by a
factor of at most 2, and therefore the convergence rate should stay
the same.

Using some additional results on random walks, Theorem \ref{th:3} leads to a tight bound (within a logarithmic factor) on the convergence time.

\begin{corollary}\label{c:1} The convergence time for the
equal-neighbor, time-invariant, symmetric model on a connected
graph on $n$ nodes, satisfies\footnote{Throughout, $\log$ will stand for the base-2 logarithm.}
\[ T_n(\epsilon) =O(n^3\log(n/\epsilon)) \]
Furthermore, for every positive integer $n$, there exists a
$n$-node connected graph for which
$$T_n(\epsilon) =\Omega(n^3 \log(1/\epsilon)).$$
\end{corollary}
\begin{proof}
The matrix $A$ is the transition probability matrix for a random
walk on the given graph, where given the current state $i$, the next
state is equally likely to be any of its neighbors (including $i$
itself). %Let $P(t)$ be the matrix $A^t$, and let
Let $p_{ij}(t)$ be the $(i,j)$th entry of the matrix $A^t$. It is
known that (Theorem 5.1\footnote{Theorem 5.1 of \cite{L96} is proved
for symmetric graphs without self-arcs. However, \alexo{the proof
does not use the absence of self-arcs, and when they are present}
the same proof yields the same result. \alexo{We refer the reader to
the derivation of Eq. (3.1) in  \cite{L96} for details}. }of
\cite{L96}),
\begin{equation} \label{timebound}
 |p_{ij}(t) - \pi_j| \leq \sqrt{\frac{d_j}{d_i}} \rho^t.
 \end{equation} \noindent
Since $1 \leq d_i$ and $ d_j\leq n$, we have
\[ |p_{ij}(t) - \pi_j| \leq \sqrt{n} \rho^t, \] for all $i$, $j$, and $t$.  Using
the result of
Theorem \ref{th:3}, we obtain
\begin{equation} |p_{ij}(t) - \pi_j| \leq \sqrt{n} (1 - n^{-3})^t.
\label{ubound}
\end{equation}
This implies that be taking $t=cn^3\log(n/\epsilon)$, where $c$ is a
sufficiently large absolute constant, we will have $|p_{ij}(\tau) -
\pi_j| \leq \epsilon/n$ for all $i$, $j$, and $\tau\geq t$.
%This implies that for $t \geq (3/4) n^3 \log n +(1/2) n^3
%\log \blue{(2/\epsilon)}$, we have $|p_{ij}(t) - \pi_j| \leq
%\epsilon/n$ for all $j$. (To see this, simply take the logarithm of
%both sides of Eq.\ (\ref{ubound}), and use the inequality $
%\log(1-n^{-3}) \geq -2n^{-3}$ for $n \geq 2$).

Let $A^* = \lim_{t \rightarrow \infty} A^t$, and
$x^*=\lim_{t\to\infty} A^t x(0)$. Note note that $A^*x(0)=x^*=A^t
x^*=A^* x^*$, for all $t$. Then, with $t$ chosen as above,
%for $t \geq (3/4) n^3 \log n + (1/2) n^3 \log (2/\epsilon)$, we have
\begin{eqnarray*}
\|x(t)-x^*\|_{\infty}&=&\|A^t(x(0)-x^*)\|_{\infty}\\
&=&\|(A^t-A^*)(x(0)-x^*)\|_{\infty}\\
&\leq& \|A^t-A^*\|_1 \cdot \|x(0)-x^*\|_{\infty}\\
&\leq& \epsilon \|x(0)-x^*\|_{\infty}.
\end{eqnarray*}
This establishes the upper bound on $T_n(\epsilon)$.

For the lower bound, note that for every $(i,j)\in\E$, we have
$\pi_i a_{ij}= (d_i/E)(1/d_i)=1/E$, so that condition (c) in Theorem \ref{th:3} is satisfied.
It follows that $A$ has real eigenvalues. Let $x(0)$ be a (real) eigenvector of $A$ corresponding to the eigenvalue $\rho$. Then, $x(t)=A^tx(0)=\rho^t x(0)$, which converges to zero, i.e., $x^*=0$. We then have
$$\frac{\|x(t)-x^*\|_{\infty}}{\| x(0)-x^* \|_{\infty}}
= \rho^t.$$
By the second part of Theorem \ref{th:3}, there exists a graph for which $\rho \geq 1-\gamma n^{-3}$, leading to the inequality $T_n(\epsilon)\geq c n^3 \log(1/\epsilon)$, for some absolute constant $c$. \end{proof}

The $\Omega(n^3)$ convergence time of this algorithm is not particularly attractive. In the next section, we explore possible improvements in the convergence time by using different choices for the matrix $A$.

\section{Convergence time for the scaled agreement algorithm\label{s:6}}

In this section, we consider the scaled agreement algorithm introduced in Section \ref{sub:scaled}.
As in \cite{XB04}, we assume the presence of a system designer who chooses the matrix $A$ so as to obtain a favorable convergence rate, subject to the condition $a_{ij}=0$ whenever $(i,j)\notin \mathcal{E}$. The latter condition is meant to represent the network topology through which the nodes are allowed to communicate. Our aim is to characterize the best possible convergence rate guarantee. We will see that the convergence rate can be brought arbitrarily close to zero. However, if we impose
a certain ``numerical stability'' requirement, the convergence time becomes $\Omega(n^2 \log (1/\epsilon))$, for a worst-case choice of the underlying graph. Furthermore, this worst-case lower bound applies even if we allow for matrices $A$ in a much larger class than that considered in \cite{XB04}.  Finally, we will show that a convergence time of $O(n^2\log (n/\epsilon))$ can be guaranteed in a simple manner, using a spanning tree.

\subsection{Favorable but impractical convergence rates} \label{span}
In this section, we show that given a connected symmetric directed
graph $G=(\N,\E)$, there is an elementary way of choosing a
stochastic matrix $A$ for which $\rho$ is arbitrarily close to zero.

We say that a  directed graph is a \emph{bidirectional spanning
tree} if (a) it is symmetric, (b) it contains all self-arcs $(i,i)$,
and (b) if we delete the self-arcs, ignore the orientation of the
arcs and remove duplicate arcs, we are left with a spanning tree.

Without loss of generality, we assume that $G$ is a bidirectional
spanning tree; since $G$ is symmetric and connected, this amounts to
deleting some of its arcs, or, equivalently, setting $a_{ij}=0$ for
all deleted arcs $(i,j)$.

Pick an arbitrary node, denoted by $r$, and designate it as the
root. Consider an arc $(i,j)$ and suppose that $j$ lies on the path
from $i$ to the root. Let $\oa_{ij}=1$ and $\oa_{ji}=0$. Finally,
let $\oa_{rr}=1$, and $\oa_{ii}=0$ for $i\neq r$. This corresponds
to a Markov chain in which the state moves deterministically towards
the root. We have ${\oA}^{t}=e_r {\bf 1}^{T}$, for all $t\geq n$,
where $e_i$ is the $i$th basis vector. It follows that $\rho=0$, and
$T_n(\epsilon)\leq n$. However, this matrix $\oA$ is not useful,
because the corresponding vector of steady-state probabilities has
mostly zero entries, which prohibits the scaling discussed in
Section \ref{sub:scaled}. Nevertheless, this is easily remedied by
perturbing the matrix $\oA$, as follows. For every $(i,j)\in \E$
with $i\neq j$ and $\oa_{ij}=0$, let $a_{ij}=\delta$, where $\delta$
is a small positive constant. For every $i$, there exists a unique
$j$ for which $\oa_{ij}=1$. For any such pair $(i,j)$, we set
$a_{ij}=1-\sum_{k=1}^n a_{ik}$ (which is nonnegative as long as
$\delta$ is chosen small enough). We have thus constructed a new
matrix $A_{\delta}$ which corresponds to a Markov chain whose
transition diagram is a
bidirectional spanning tree. %, together with self-arcs.
Since \alexo{the convergence rate $\rho$ is an eigenvalue of the
iteration matrix, and }eigenvalues are continuous functions of
matrix elements, we see that, for the matrix $A_{\delta}$, the
convergence rate $\rho$ can be made as small as desired, by choosing
$\delta$ sufficiently small. Finally, since $A_{\delta}$ is a
positive matrix,  the corresponding vector of steady-state
probabilities is positive.

To summarize, by picking $\delta$ suitably small, we can choose a
(stochastic) matrix $A_{\delta}$ with arbitrarily favorable
convergence rate, and which allows the application of the scaled
agreement  algorithm of Section \ref{sub:scaled}. It can be shown
that the convergence time is linear in the following sense: For
every $\epsilon$, there exists some $\delta$ such that, for the
matrix $A_{\delta}$, the corresponding convergence time, denoted by
$T_n(\epsilon; \delta)$, satisfies $T_n(\epsilon; \delta) \leq n$.
Indeed, this is an easy consequence of the facts $\lim_{\delta\to
0}(A_{\delta}^n -\oA^n)=0$ and $T_n(\epsilon';0)\leq n$ for every
$\epsilon'>0$\footnote{\alexo{Indeed, it is easy to see that by
suitably choosing the root, we can make sure that convergence time
is at most $\lceil d(G)/2 \rceil$ where $d(G)$ is the diameter of
the graph $G$ defined as the largest distance between any two
vertices.} }.

%Indeed, let us denote by $P$ the limiting
%chain and by $P_{\delta}$ the approximate chain parametrized by
%$\delta$. It can be shown that $\lim_{\delta \rightarrow 0}
%P_{\delta}^{n} - P_{\delta}^{\infty} = P^n - P^{\infty}=0$,
%suggesting that we may make the error at time $n$ arbitrarily small
%by picking a small enough $\delta$.

However, note that as $n$ gets larger, $n \pi_i$ may approach $0$ at
the non-root nodes. The implementation of the scaling in Eq.\
(\ref{eq:scale}) will involve division by a number which approaches
$0$, possibly leading to numerical difficulties. Thus, the resulting
averaging algorithm may be undesirable. Setting averaging aside, the
agreement algorithm based on $A_{\delta}$, with $\delta$ small is
also undesirable: despite its favorable convergence rate, the final
value on which consensus is reached is approximately equal to the
initial value $x_r(0)$ of the root node. Such a ``dictatorial''
solution runs contrary to the motivation behind consensus
algorithms.
% in numerical implementation for large $n$.

\subsection{A lower bound}

In order to avoid the numerical issues raised above,
we will now impose a condition on the dominant (and positive) left eigenvector $\pi$ of the matrix $A$, and require
\begin{equation} \label{cbound}
{n \pi_i \geq \frac{1}{C}, \quad\forall ~i}\end{equation} where $C$
is a given constant with $C > 1$. This condition ensures that $n
\pi_i$ does not approach $0$ as $n$ gets large, so that the initial
conditions in the scaled agreement algorithm of Section
\ref{sub:scaled} are well-behaved. Furthermore, in the context of
consensus algorithms, condition (\ref{cbound}) has an appealing
interpretation: it requires that the initial value $x_i(0)$ of every
node $i$ has a nonnegligible impact on the final value
$\lim_{t\to\infty} x_k(t)$, on which consensus is
reached\footnote{\alexo{In the case where $A$ is the transition
matrix of a reversible Markov chain, there is an additional
interpretation. A reversible Markov chain may be viewed as a random
walk on an undirected graph with edge-weights. Defining the degree
of an vertex as the sum total of of the weights incident upon it,
the condition $n \pi_i \geq C$ is equivalent to requiring that each
degree is lower bounded by a constant times the average degree.}}.

We will now show that under the additional condition (\ref{cbound}),
there are graphs for which the convergence time is $\Omega(n^2
\log(1/\epsilon))$. One may wonder whether a better convergence time
is possible by allowing some of the entries of $A$ to be negative.
As the following result shows, negative entries do not help. The
graph that we employ is a {\em line graph}, with arc  set
$\E=\{(i,j)\mid |i-j|\leq 1\}$.

\begin{theorem}\label{th:4}
Consider an $n\times n$ matrix $A$ such that $a_{ij}=0$ whenever
$|i-j| > 1$, and such that the graph with edge set
$\{(i,j) \in \E ~|~ a_{ij} \neq 0\}$ is connected. Let
$\lambda_1,\lambda_2,\ldots$, be its eigenvalues in order of
decreasing modulus. Suppose that $\lambda_1=1$ and $A{\bf 1}=1$.
Furthermore, suppose that there exists a vector $\pi$ satisfying
Eq.\ (\ref{cbound}) such that $\pi^T A=\pi^T$. Then, there exist
absolute constants $c_1$ and $c_2$ such that
\[ \rho \geq 1 - c_1\frac{C}{n^2}, \]
and
$$T_n(\epsilon) \geq c_2\frac{n^2}{C} \log({1}/{\epsilon}).$$
\end{theorem}\begin{proof}
If the entries of $A$ were all nonnegative, we would be dealing with a birth-death Markov chain. Such a chain is reversible, i.e., satisfies the detailed balance equations $\pi_ia_{ij}= \pi_j a_{ji}$ (condition (c) in Theorem \ref{th:2}). In fact the derivation of the detailed balance equations does not make use of nonnegativity; thus, detailed balance holds in our case as well.

Without loss of generality, we can assume that $\sin \pi_i =1$.
For $i=1,\ldots,n$, let $y_i=i-\beta$, where $\beta$ is chosen so that
$\sin \pi_i y_i=0$. We will make use of the inequality (\ref{mineigensym}).
Since $a_{ij}=0$ whenever $|i-j|>1$, we have
\begin{equation}\sin\sjn \pi_i a_{ij} (y_i-y_j)^2 \leq \sin\sjn \pi_i a_{ij} =1.\label{eq:aa}\end{equation}
Furthermore,
\begin{equation}\label{eq:bb}
\sin \pi_i y_i^2
\geq
\frac{1}{nC} \sin y_i^2
=
 \frac{1}{nC} \sin (i-\beta)^2
 \geq
  \frac{1}{nC} \sin \Big(i- \frac{n+1}{2}\Big)^2
  \geq
  \frac{n^2}{12 C}.\end{equation}
The next to last inequality above is an instance of the general
inequality ${\mathbf E}[(X-\beta)^2] \geq {\rm var}(X)$, applied to
a discrete uniform random variable $X$. The last inequality follows
from the well known fact ${\rm var}(X)= (n^2-1)/12$.
%Let us now bound $\beta$. We have $\beta=\sin \pi_i i$, and $\pi$ is
%constrained by \red{$\pi_i \geq 1/nC$. Thus, \[ \beta \geq
%\frac{n(n+1)}{2} \frac{1}{nC} >  \frac{n}{2C} \] For every $i\leq
%n/4C$, we have $y_i^2=(\beta-i)^2 >(n/4 C)^2$. Since $ \pi_i \geq
%1/Cn$, we conclude that
%\begin{equation}
%\sin \pi_i y_i^2 \geq \sum_{i\leq n/4 C} \pi_i y_i^2 \geq \frac{n}{4
%C}\cdot \frac{1}{nC} \cdot \frac{n^2}{4^2 C^2} = \frac{n^2}{4^3
%C^4}. \label{eq:bb}
%\end{equation}}
Using the inequality (\ref{mineigensym}) and Eqs.\ (\ref{eq:aa})-(\ref{eq:bb}), we obtain
the desired bound on $\rho$.

For the bound on $T_n(\epsilon)$, we let $x(0)$ be a (real)
eigenvector of $A$, associated with the eigenvalue $\lambda_2$, and
proceed as in the end of the proof of Corollary \ref{c:1}.
\end{proof}

\bigskip

\alexo{\noindent {\bf Remark:} Note that if the matrix $A$ is as in
the previous theorem, it is possible for the iteration
$x(t+1)=Ax(t)$ not to converge at all. Indeed, nothing in the
argument precludes the possibility that the smallest eigenvalue is
$-1$, for example. In such a case, the lower bounds of the theorem -
derived based on bounding the second largest eigenvalue - still hold
as the convergence rate and time are infinite.}

\subsection{Convergence time for spanning trees\label{treeher}}
We finally show that an $O(n^2)$ convergence time guarantee is easily obtained, by restricting to a spanning tree.
\begin{theorem}\label{th:5}
Consider the equal-neighbor, time-invariant, bidirectional model on
a bidirectional spanning tree. We have
\[ \rho \leq 1 - \frac{1}{3 n^2}, \]
and
\[ T_{n} (\epsilon) =O\big(n^2 \log({n}/{\epsilon})\big).\]
\end{theorem}
\begin{proof}
In this context, we have $\pi_i=d_i/E$, where
$E=\sin d_i =2(n-1)+n<3n$. (The factor of 2 is because we have arcs in both directions; the additional term $n$ corresponds to the self-arcs.)
As in the proof of Theorem \ref{th:4}, the detailed balance conditions $\pi a_{ij}=\pi_j a_{ji}$ hold, and we can apply Theorem \ref{th:2}.
Note that Eq.\ (\ref{eq:gg})
can be rewritten in the form
\begin{equation}
\label{noneigensym} \lambda_2 = 1 - \frac{1}{2} \min_{\sum_i^n d_i x_i
= 0, \sum_i^n d_i x_i^2 = 1} \sum_{(i,j)\in\E} (x_i - x_j)^2.
\end{equation}

We use the methods of \cite{LO81} to show that for trees,
$\lambda_2$ can be upper bounded by $1 - 1/3n^2$. Indeed, suppose
that $x$ satisfies $\sum_i^n d_i x_i = 0$ and $\sum_i^n d_i x_i^2 = 1$,
and let $x_{\rm \max}$ be such that $|x_{\rm max}|=\max_i |x_i|$.
Then,
\[ 1 = \sum_i d_i x_i^2 \leq  3n x_{\rm max}^2, \] and
it follows that $|x_{\max}| \geq 1/\sqrt{3n}$. Without loss of
generality, assume that $x_{\rm max} > 0$ (else, replace each $x_i$ by
$-x_i$).  Since $\sum_i d_i x_i = 0$, there exists some $i$ for
which $x_i<0$; let us denote such a negative $x_i$ by $x_{\rm neg}$.
Then,
\begin{equation}  \label{chain} \frac{1}{\sqrt{3n}} \leq x_{\max} - x_{\rm neg}  =
(x_{\max} - x_{k_1}) + (x_{k_1} - x_{k_2}) + \cdots   + (x_{k_{r-1}}
- x_{\rm neg}),
\end{equation}
where $k_1,k_2,\ldots,k_{r-1}$ are the nodes on the path from
$x_{\max}$ to $x_{\rm neg}$. By the Cauchy-Schwartz inequality,
\begin{equation} \label{diffbound} \frac{1}{3n} \leq \frac{n}{2} \sum_{(i,j) \in
\mathcal{E}} (x_i - x_j)^2.
\end{equation}
(The
factor of 1/2 in the right-hand side arises because the sum includes
both terms $(x_{k_i}-x_{k_{i+1}})^2$ and
$(x_{k_{i+1}}-x_{k_{i}})^2$.)  Thus,
\[\sum_{(i,j) \in \mathcal{E}} (x_i - x_j)^2 \geq \frac{2}{3n^2}. \]
which proves the bound for the second largest  eigenvalue.

For the smallest eigenvalue, recall that $a_{ii}\geq 1/n$ for every
$i$. It follows that the matrix $A$ is of the form $I/n + Q$, where
$I$ is the identity matrix, and $Q$ is a nonnegative matrix whose
row sums are equal to $1-1/n$.  Thus, all of the eigenvalues of $Q$
have magnitude bounded above by $1-1/n$, which implies that the
smallest eigenvalue of $Q$ is bounded below by $-1+1/n$. We conclude
that $\lambda_n$, the smallest eigenvalue of $I/n+Q$, satisfies
$$\lambda_n \geq -1 +\frac{2}{n} \geq -1 +\frac{2}{n^3}.$$

For the bound on the convergence time, we proceed as in the proof of
Corollary \ref{c:1}. Let $p_{ij}(t)$ be the $(i,j)$th entry of
$A^t$. Then,
\[ |p_{ij}(t) - \pi_j| \leq \sqrt{n} \Big(1 - \frac{1}{3} n^{-2}\Big)^t. \]
For a suitable absolute constant $c$ and for $t \geq c n^2 \log
(n/\epsilon)$, we obtain  $|p_{ij}(t) - \pi(j)| \leq
{\epsilon}/{n}$. The rest of the proof of Corollary \ref{c:1} goes
through unchanged. \end{proof}

In light of the preceding theorem, we propose the following simple heuristic, with worst-case convergence time $O(n^2\log(n/\epsilon))$, as an alternative to a more elaborate design of the matrix $A$.

\begin{algorithm}\label{al:2} {\rm We are given a symmetric graph $G$. We
delete enough arcs to turn $G$ into a bidirectional spanning tree,
and then carry out the} equal-neighbor, time-invariant,
bidirectional consensus algorithm, with initial value
$x_i(0)/n\pi_i$ at node~$i$.
\end{algorithm}

Let us remark that the $O(n^2\log(n/\epsilon))$ bound (Theorem
\ref{th:5}) on the convergence time of this heuristic is essentially
tight (within a factor of $\log n$). Indeed, if the given graph is a
line graph, then with our heuristic we have $n \pi_i = n d_i/E \geq
2/3$, and Theorem \ref{th:4} provides a
$\Omega(n^2\log(1/\epsilon))$ lower bound.

\section{Convergence time when using a doubly stochastic matrix} \label{s:7}

We provide here a brief comparison of our methods with two methods that have been proposed in the literature, and which rely on doubly stochastic matrices. Recall that doubly stochastic matrices give rise directly to an averaging algorithm, without the need for scaling the initial values.

\begin{itemize}
\item[(a)]
Reference \cite{XB04} considers the case where the graph $G$ is given, and studies the problem of choosing a doubly stochastic matrix
$A$ for which the convergence rate $\rho$ is smallest. In order to obtain a tractable (semidefinite programming) formulation, this reference imposes the further restriction that $A$ be symmetric. For a doubly stochastic matrix, we have $\pi_i=1/n$, for all $i$, so that condition (\ref{cbound}) holds with $C=1$. According to Theorem \ref{th:4}, there exists a sequence of graphs, for which we have $T_n(\epsilon)=\Omega(n^2\log (1/\epsilon))$. We conclude that despite the sophistication of this method, its worst case guarantee is no better (ignoring the $\log n$ factor) than the simple heuristic we have proposed (Algorithm \ref{al:2}). On the other hand, for particular graphs, the design method of \cite{XB04} may yield better convergence times.
\item[(b)]
The following method
was proposed in \cite{OSM04}. The nodes first agree on some
value $\epsilon \in (0,1/\max_i d_i)$. (This is easily accomplished in a distributed manner.)
%example, computing the  largest degree $d_{\rm max}$ and agreeing on
%$1/(2 d_{\rm max})$. The largest degree could be computed by a
%simple iterative algorithm: at each step, each node takes the
%maximum of its own degree and that of each of its neighbors. After
%at most $n$ steps, each node knows the highest degree.
Then, the nodes iterate according to
the equation
\begin{equation} \label{osmiter}
x_i(t+1) = (1 - \epsilon d_i) x_i(t) + \epsilon \sum_{j \in \N(i) \setminus
\{i\}}^n x_j(t). \end{equation}
Assuming a connected graph, the iteration converges
to
consensus (this is a special case of Theorem \ref{th:1}).
Furthermore, this iteration preserves the sum $\sin x_i(t)$. Equivalently,
the corresponding matrix $A$ is doubly stochastic, as required in order to have an averaging algorithm.

This algorithm has the disadvantage of uniformly small step sizes.
If many of the nodes have degrees of the order of $n$, there is no
significant theoretical difference between this approach and our
Algorithm \ref{al:1}, as both have effective step sizes of order of
$1/n$. On the other hand, if only a small number of nodes have large
degree, then the algorithm in \cite{OSM04} will force {\em all} the
nodes to take small steps. This drawback is avoided by our
Algorithms \ref{al:1} (Section  \ref{3a}) and \ref{al:2} (Section
\ref{treeher}). A comparison of the method of \cite{OSM04} with
Algorithm \ref{al:1} is carried out, through simulation experiments,
in Section \ref{s:8}.
\end{itemize}

\section{Averaging with dynamic topologies\label{lb}}\label{s:8}

In this section, we turn our attention to the more challenging case where communications are bidirectional but the network topology changes dynamically. Averaging algorithms for such a context have been considered
previously in
\cite{MSPLM05,MVR05}.

As should be clear from the previous sections, consensus and
averaging algorithms are intimately linked, with the agreement
algorithm often providing a foundation for the development of an
averaging algorithm. For this reason, we start by investigating the
worst-case performance of the agreement algorithm in a dynamic
environment. Unfortunately, as shown in Section \ref{expon}, its
convergence time is not polynomially bounded, in general, even
though it is an open question whether this is also the case when we
restrict to symmetric graphs. Motivated by this negative result, we
approach the averaging problem differently: we introduce an
averaging algorithm based on ``load balancing'' ideas (Section
\ref{desc}), and prove a polynomial bound on its convergence time
(Section \ref{con}).

\subsection{Non-polynomial convergence time for the agreement algorithm\label{expon}}

We begin by formally defining the notion of ``convergence time'' for
the agreement algorithm on dynamic graph sequences. Given a sequence
of graphs $G(t)$ on $n$ nodes such that Assumption \ref{as:3} of
Section \ref{s:2} is satisfied for some $B > 0$, and an initial
condition $x(0)$,  we define the convergence time
$T_{G(\cdot)}(x(0),\epsilon)$ (for this particular graph sequence
and initial condition) as the first time $t$ when each node is
within \alexo{an $\epsilon$-neighborhood} of the final consensus,
i.e., $\|x(t) - \lim_{t\to\infty} x(t) \|_{\infty} \leq \epsilon$.
We then define the (worst-case) convergence time, $T_n(B,\epsilon)$,
as the maximum value of $T_{G(\cdot)}(x(0),\epsilon)$, over all
graph sequences $G(\cdot)$ on $n$ nodes that satisfy Assumption
\ref{as:3} for that particular $B$, and all initial conditions that
satisfy $\|x(0)\|_{\infty}\leq 1$.

%Thus the
%convergence time depends not only on the tolerance $\epsilon$ and
%the initial distribution of values $x(0)$, but also on the graph
%sequence $G(\cdot)$.

We focus our attention on the equal-neighbor version of the
agreement algorithm. The result that follows shows that its
convergence time is not bounded by a polynomial in $n$ and $B$. In
particular, if $B$ is proportional to $n$, the convergence time
increases faster than an exponential in $n$. We note that the upper
bound in Theorem \ref{th:dyn} is not a new result, but we include it
for completeness, and for comparison with the lower bound, together
with a proof sketch. Similar upper bounds have also been provided
recently in \cite{CSM05}, under slightly different assumptions
on the graph sequence $G(\cdot)$.

\begin{theorem}\label{th:dyn}
For the equal-neighbor agreement algorithm, there exist positive
constants $c_1$ and $c_2$ such that for every $n$, $B$, and
$\alexo{1>}\epsilon>0$,
\begin{equation}\label{eq:dynth}
 c_1 nB\Big(\frac{n-1}{2}\Big)^{B-1} \log \frac{1}{\epsilon}
  \leq T_n(B,\epsilon) \leq c_2 B n^{nB} \log \frac{1}{\epsilon}.
  \end{equation}
  \end{theorem}
\begin{proof}
The upper bound follows by inspecting the proof of convergence of
the agreement algorithm with the constant $\alpha$ of Assumption
\ref{as:1} set to $1/n$ (c.f. \cite{T84, BHOT05}).

\old{We will show that there exists a sequence ${G}(t)$ of directed
graphs with $B=n+1$, and a corresponding set of initial values
${x}(0)$, such that $T_{{G}(t)}({x}(0), \epsilon)=\Omega((n/2)^{n}
\log ({1}/{\epsilon}))$. On the other hand, by inspecting the proof
of convergence of the agreement algorithm, it is easy to see that
the convergence time admits a $O(n^{n} B \log ({1}/{\epsilon}))$
upper bound
--- see \cite{CSM05} and \cite{BHOT05}. We can
therefore conclude that these upper bounds are essentially tight.}

We now prove the lower bound by exhibiting a sequence of graphs
$G(t)$ and an initial vector $x(0)$, with $\|x(0)\|_{\infty}\leq 1$
for which $T_{G(\cdot)} (x(0),\epsilon) \geq c_1 nB
(n/2)^{B-1}\log(1/\epsilon)$. We assume that $n$ is even and $n\geq
4$. The initial condition ${x}(0)$ is defined as ${x}_i(0)=1$ for $i
= 1,\ldots,n/2$, and ${x}_i(0) = -1$ for $i=n/2+1,\ldots,n$.
\begin{itemize}
\item[(i)] The graph $G(0)$, used for the first iteration, is shown in the left-hand side of Figure \ref{f:diag}.

\begin{figure}\label{f:diag}
\includegraphics[width=6cm]{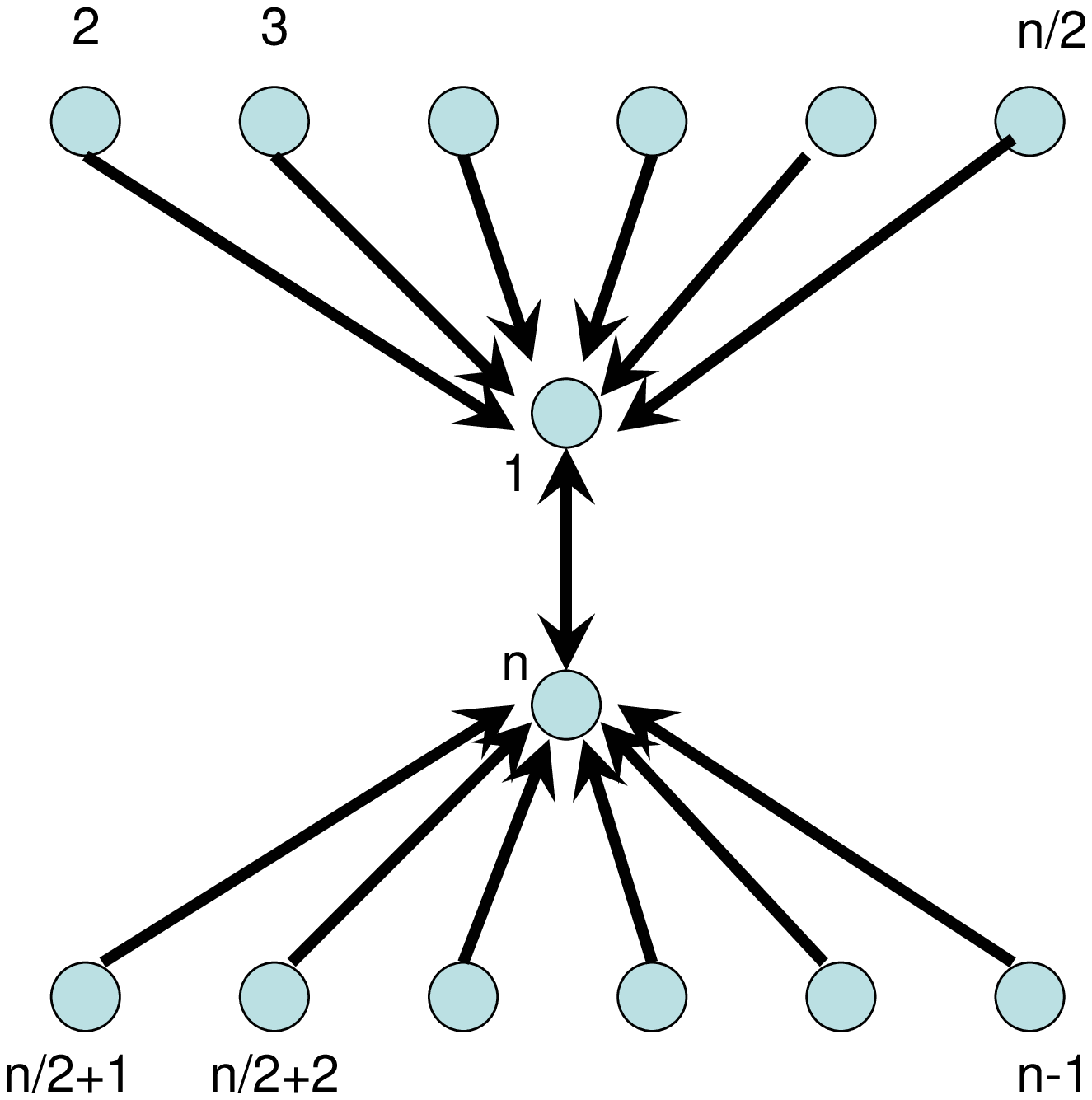}\hspace{1cm}\includegraphics[width=6cm]{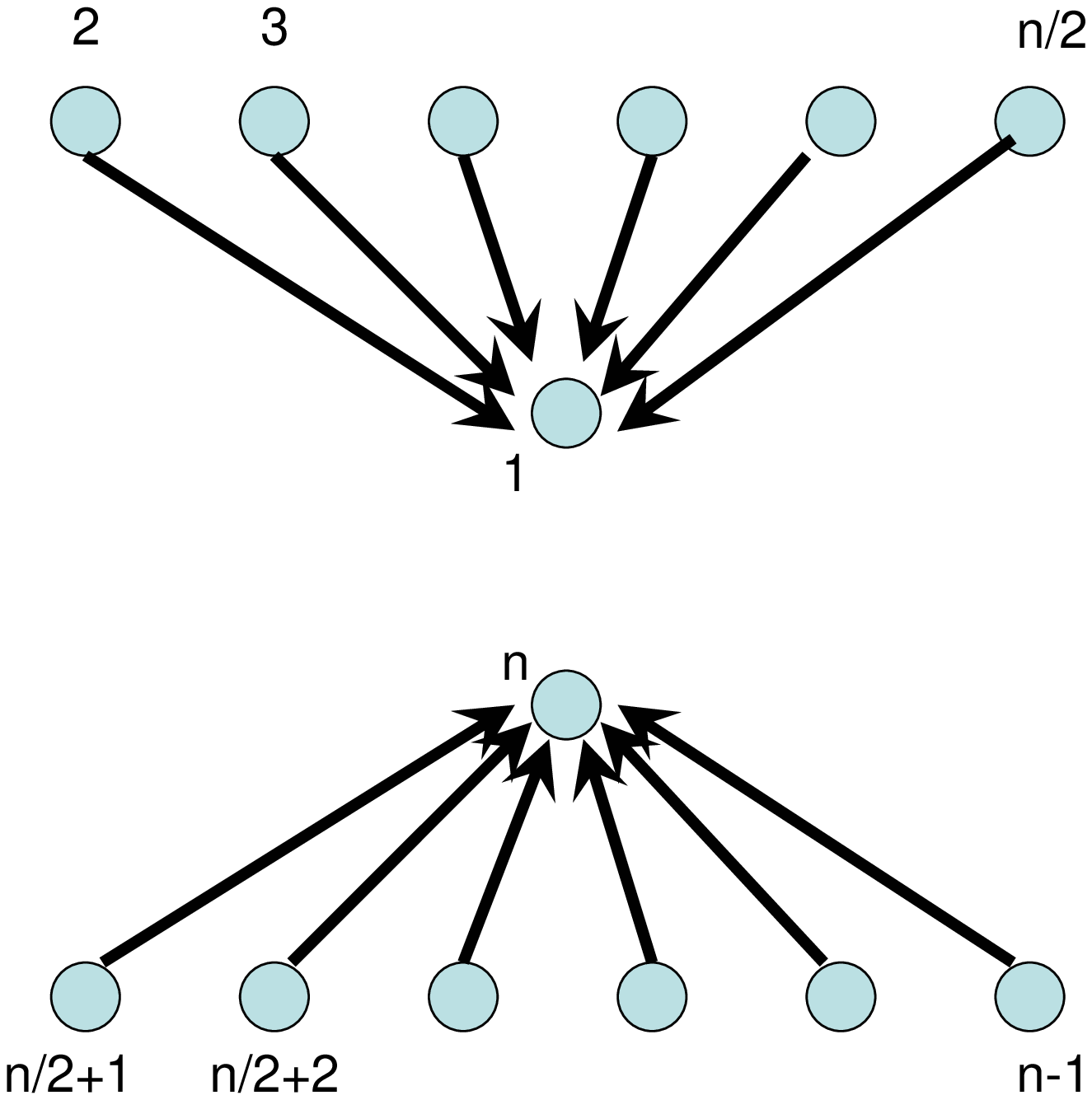}
\caption{\label{t0} The diagram on the left is the graph $G(0)$. The
diagram on the right is the graph $G(t)$ at times $t = 1,\ldots,
B-2$. Self-arcs are not drawn but should be assumed present at every
node.}
\end{figure}

\item[(ii)] For $t = 1,\ldots, B-2$, we perform an equal-neighbor iteration, each time using the graph $G(t)$ shown in the right-hand side of
Figure \ref{f:diag}.

\item[(iii)] Finally, at time $B-1$,
the graph $G(B-1)$ consists of the complete graph over the
nodes $\{1,\ldots,n/2\}$ and the complete graph over the nodes
$\{n/2+1,\ldots,n\}$.
\item[(iv)] This sequence of $B$ graphs is then repeated: $G(t+kB)=G(t)$ for every positive integer $k$.
\end{itemize}
It is easily seen that this sequence of graphs satisfies Assumption \ref{as:3}, and that convergence to consensus is guaranteed.

At the end of the first iteration, we have $x_i(1)=x_i(0)$, for $i\neq 1,n$,
and
\begin{equation}\label{eq:alphao}
x_1(1)=\frac{(n/2)-1}{(n/2)+1} =1- \frac{4}{n+2},\qquad x_n(1)
=-x_1(1).
\end{equation}
Consider now the evolution of $x_1(t)$, for $t=1,\ldots,B-2$, and
let $\alpha(t)=1-x_1(t)$. We have
\[ x_1(t+1) = \frac{1 \cdot (1-\alpha(t)) + (n/2-1) \cdot
1}{n/2} = 1 - (2/n) \alpha(t), \] so that
$\alpha(t+1)=2\alpha(t)/n$. From Eq.\ (\ref{eq:alphao}),
$\alpha(1)=4/(n+2)$, which implies that $\alpha(B-1) =(2/n)^{B-2}$,
or
\[ x_1(B-1) = 1-\frac{4}{n+2} \Big(\frac{2}{n}\Big)^{B-2} \] By symmetry, \[
x_n(B-1) = -1 + \frac{4}{n+2} \Big(\frac{2}{n}\Big)^{B-2}. \]

Finally, at time $B-1$, we iterate on the complete graph over nodes
$\{1,\ldots,n/2\}$ and the complete graph over nodes
$\{n/2+1,\ldots,n\}$. For $i=2,\ldots,n/2$, we have $x_i(B-1)=1$,
and  we obtain
  \[x_i(B-1)=
\frac{
\displaystyle{1\cdot\Big(\frac{n}{2} -1\Big) + 1-\frac{4}{n+2} \Big(\frac{2}{n}\Big)^{B-2}}}
 {n/2} =
1-\frac{4}{n+2} \Big(\frac{2}{n}\Big)^{B-1}.\] Similarly, for
$i=(n/2)+1,\ldots,n$, we obtain  \[ x_i(B-1)= -1 + \frac{4}{n+2}
\Big(\frac{2}{n}\Big)^{B-2}.\]

Thus,
\[ \frac{|\max_i {x}_i(B) - \min_i {x}_i(B)|}{|\max_i
{x}_i(0) - \min_i {x}_i(0)|}=  1 - \frac{4}{n+2} \cdot \Big(\frac{2}{n}\Big)^{B-1}.
\]
Moreover, because ${x}(B)$ is simply a scaled version of ${x}(0)$, it
is clear that by repeating this sequence of graphs, we will have
$$\frac{|\max_i
{x}_i(kB) - \min_i {x}_i(kB)|}{|\max_i {x}_i(0) - \min_i
{x}_i(0)|}=\Big( 1 - \frac{4}{n+2} \cdot \Big(\frac{2}{n}\Big)^{B-1}
\Big)^{k}.$$ This readily implies that
\[ T_{G(\cdot)(t)}({x}(0),\epsilon) =\Omega\Big(nB\Big(\frac{n}{2}\Big)^{B-1} \log \frac{1}{\epsilon}
\Big). \]

If $n$ is odd, then $n' = n-1$ is even. We apply the
same initial condition and graph sequence as above to nodes $\{1,
\ldots, n'\}$. As for the additional node $x_n$, we let $x_n(0)=0$
and make extra connections by connecting node $n$ to nodes $1$ and
$n'$ at time $0$ with a bidirectional link. By repeating the
analysis above, it can be verified that
\[ T_{G(\cdot)(t)}({x}(0),\epsilon)
=\Omega\Big(nB\Big(\frac{n-1}{2}\Big)^{B-1} \log \frac{1}{\epsilon}
\Big). \] This concludes the proof.
\old{In particular, the
convergence time cannot be bounded by a polynomial function of $n$
and $B$. Furthermore, if $B$ is proportional to $n$, the convergence
time becomes $\Omega\big((n/2)^n\big)$.}
\end{proof}

Both the upper and lower bounds in Theorem \ref{th:dyn} display an
exponential growth of the convergence time, as a function of $B$. It
is unclear, however, which of the two terms, $n^B$ or $n^{nB}$,
better captures the behavior of $T_n(B,\epsilon)$.

\subsection{Polynomial-time averaging in dynamic topologies\label{desc}}

The algorithm we present here is a variation of an old {\em load
balancing} algorithm (see \cite{C89} and Chapter 7.3 of
\cite{BT89}). Intuitively, a collection of processors with different
initial loads try to equalize their respective loads. As some of the
highly loaded processors send some of their load to their less
loaded neighbors, the loads at different nodes tend to become equal.
Similarly, at each step of our algorithm, each node offers some of
its value to its neighbors, and accepts or rejects such offers from
its neighbors. Once an offer from $i$ to $j$ to send $\delta$ has
been accepted, the updates $x_i := x_i - \delta$ and $x_j := x_j +
\delta$ are executed.

We assume a time-varying sequence of graphs $G(t)$. We only make two
assumptions on $G(\cdot)$: symmetry and bounded interconnectivity
times (see Section \ref{s:2} for definitions). The symmetry
assumption is natural if we consider, for example, communication
between two nodes to be feasible whenever the nodes are within a
certain distance of each other. The assumption of bounded
interconnectivity times is necessary for an upper bound on the
convergence time to exist (otherwise, we could insert infinitely
many empty graphs $G(t)$, in which case convergence is arbitrarily
slow for any algorithm).

We next describe formally the steps that each node carries out at
each time $t$. For definiteness, we refer to the node executing the
steps below as node $A$. Moreover, the instructions below sometimes
refer to the ``neighbors'' of node $A$; this always means nodes
other than $A$ that are neighbors at time $t$, when the step is
being executed (since $G(t)$ can change with $t$,  the set of
neighbors of $A$ can also change). Let $\N_i(t)=\{j\neq i: (i,j)\in
\E(t)\}$. Note that this a little different from the definition of
$\N_i(t)$ in earlier sections, in that $i$ is no longer considered a
neighbor of itself.
\begin{algorithm}\label{al:3}
\rm
If $\N_A(t)$ is empty, node $A$ does nothing at time $t$. Else, node $A$ carries out the following steps.
\begin{enumerate}
\item[1.]  Node $A$ broadcasts its current value $x_A$ to all of its
neighbors (every $j$ with $j\in \N_A(t)$).

\item[2.] Node $A$ finds a neighboring node $B$ with the smallest value: $x_B=\min\{x_j: j\in\N_A(t)\}$. If $ x_A\leq x_B$, then node $A$ does nothing further at
this step. If $x_B< x_A$, then
node $A$ makes an offer of
$ (x_A - x_B)/2$ to node $B$.

\item[3.] If node $A$ does not receive any offers, it does nothing further at this step. Otherwise,
it sends an
acceptance to the sender of the largest offer and a rejection to all
the other senders. It updates the value of $x_A$ by adding the value
of the accepted offer.

\item[4.] If an acceptance arrives for the offer made by node
$A$, node $A$ updates $x_A$ by subtracting the value of the offer.
\end{enumerate}
\end{algorithm}

For concreteness, we use $x_i(t)$ to denote the value possessed by
node $i$ at the {\em beginning} of the above described steps.
Accordingly, the value possessed by node $i$ at the end of the above
steps will be $x_i(t+1)$. The algorithm we have specified clearly
keeps the value of $\sin x_i(t)$ constant. Furthermore, it is a
valid averaging algorithm, as stated in Theorem \ref{th:6} below. We
do not provide a separate proof, because this result follows from
the convergence time bounds in the next subsection.
\begin{theorem}\label{th:6}
Suppose that each $G(t)$ is symmetric and that Assumption \ref{as:3}
(bounded interconnectivity times) holds. Then, $\lim_{t \rightarrow
\infty} x_i(t) = \frac{1}{n} \sum_{k=1}^n x_k(0)$, for all~$i$.
\end{theorem}

\subsection{Convergence time\label{con}}
We introduce the following ``Lyapunov'' function that quantifies the distance of the state $x(t)$ of the agents form the desired limit:
$$V(t)= \Big\|{x}(t) - \frac{1}{n} \sum_{i=1}^n
x_i(0) {\bf 1}\Big\|_2^2.$$ Intuitively, $V(t)$ measures the
variance of the values at the different nodes. Given a sequence of
graphs $G(t)$ on $n$ nodes, and an initial vector $x(0)$, we define
the convergence time $T_{G(\cdot)}(x(0),\epsilon)$ as the first time
$t$ after which $V(\cdot)$ remains smaller than $\epsilon V(0)$:
$$T_{G(\cdot)}(x(0),\epsilon)=\min\big\{t\mid V(\tau)\leq \epsilon V(0),\ \forall\ \tau\geq t\big\}.$$
We then define the (worst-case) convergence time, $T_n(B,\epsilon)$,
as the maximum value of $T_{G(\cdot)}(x(0),\epsilon)$, over all
graph sequences $G(\cdot)$ on $n$ nodes that satisfy Assumption
\ref{as:3} for that particular $B$, and all initial conditions
$x(0)$.

\begin{theorem} There exists
a constant $c>0$ such that for every $n$ and $\alexo{1>}\epsilon>0$,  %and every
%sequence of symmetric \red{graphs} $G(t)$ that satisfy Assumption
\ref{as:3}, we have
\begin{equation}
T_n(B, \epsilon) \leq c B n^3 \log \frac{1}{\epsilon}.\label{ncubed}
\end{equation}

\end{theorem}
\begin{proof}
The proof is structured as follows.
Without loss of generality, we assume that $\sin
x_i(0) = 0$; this is possible because adding a constant to each $x_i$ does not change the sizes of the offers or the acceptance decisions. We will show that $V(t)$ is nonincreasing in
$t$, and that
\begin{equation} \label{eq:vdif} { V}((k+1)B)  \leq  \Big(1 - \frac{1}{2 n^3}\Big){ V}(kB) \end{equation}
for every nonnegative integer $k$. These two claims readily imply
the desired result. To see this, note that if $V(t)$ decreases by a
factor of $1-({1}/{2 n^3})$ every $B$ steps, then it decreases by a
$\Theta(1)$ factor in $Bn^3$ steps. It follows that the time until
$V(t)$ becomes less than $\epsilon V(0)$ is $O(Bn^3
\log(1/\epsilon))$. Finally, since $V(t)$ is nonincreasing, $V(t)$
stays below $\epsilon V(0)$ thereafter.

We first show that $V(t)$ is nonincreasing. We argue that while
rejected offers clearly do not change $V(t)$, each accepted offer at
time $t$ results in a decrease of $V(t+1)$. While this would
straightforward to establish if there were a single accepted offer,
a more complicated argument is needed to account for the possibility
of multiple offers being simultaneously accepted. We will show that
we can view the changes at time $t$ as a result of a series of
sequentially accepted offers, each of which results in a smaller
value of $V$.

Let us focus on a particular time $t$. We order the nodes from
smallest to largest, so that $x_1(t) \leq x_2(t) \leq \cdots \leq
x_n(t)$, breaking ties arbitrarily. Let $A_i(t)$ be the size of the
offer accepted by node $i$ at time $t$ (if any). If the node
accepted no offers at time $t$, set $A_i(t)=0$. Furthermore, if
$A_i(t)> 0$, let $\mathcal{A}_i(t)$ be the the index of the node
whose offer node $i$ accepted.

Let us now break time $t$ into $n$ %separate
{\em periods}. The $i$th
period involves the updates caused by node $i$ accepting an offer
from node $\A_i(t)$. In particular, node ${i}$  performs the
update $x_{{i}}(t):=  x_{{i}}(t) + A_{{i}}(t)$ and node
$\mathcal{A}_{{i}}(t)$ performs the update
$x_{\mathcal{A}_{{i}}(t)} (t) := x_{\mathcal{A}_{{i}}(t)}
(t) - A_{{i}}(t)$.

We note that every offer accepted at time $t$ appears
in some %time
period in the above sequence. We next argue that each
offer decreases $V$. This will complete the proof that $V(t)$ is
nonincreasing in $t$.

Let us suppose that in the $i$th period, node $i$ accepts an offer
from node $\mathcal{A}_{{i}}(t)$, \alexo{which we will for
simplicity denote by $j$}.
%after breaking down time $t$ into
%$n$ time periods as described above, an offer occurs at the $i$th
%period. Then the offer from some node must be accepted by node $i$;
%let us refer to the sender of the offer as node $j$.
%Now,
Because nodes only send offers
to lower valued nodes, the inequality $x_j > x_i$ must hold at the beginning of time
$t$, before time period $1$. We claim that this inequality continues to hold
when the $i$th time period is reached. Indeed, $x_j$
is unchanged during periods $1,\ldots,i-1$ (it can only send one
offer, which was to $x_i$; and if it receives %has received
any %other
offers, their effects will occur in period $j$, which is after period $i$).
Moreover, while the value of $x_i$ may have changed in periods $1,\ldots,i-1$,
it cannot have increased (since $i$ is not allowed to accept more than one offer at any given time $t$).
Therefore, the inequality $x_j>x_i$ still holds at the beginning of the $i$th period.

During the $i$th period, a certain positive amount is transferred from node $j$ to node $i$. Since the transfer takes place from a higher-valued node to a lower-valued one, it is easily checked that the value of
$x_i^2+x_j^2$ (which is the contribution of these two nodes to $V$) is reduced.
To summarize, we have shown that we can serialize the offers accepted
at time $t$, in such a way that each accepted offer causes a reduction in
$V$. It follows that $V(t)$ is nonincreasing.

We will now argue that at some time $t$ in the interval $0,1,\ldots,B-1$, there will be some update (acceptance of an offer) that reduces $V(t)$  by
at least $1/(2 n^3)V(0)$. Without loss of generality, we assume
$ \max_i |x_i(0)|=1$, so that all the values lie in the interval
$[-1,+1]$. It follows that $V(0) \leq n$.

%Let us order the nodes in increasing  order of values at time $0$ so
%that $ x_1(0) \leq x_2(0)  \leq \ldots \leq x_n(0)$.
Since $\sin x_i(0) = 0$, it follows that $\min_i x_i(0) \leq 0$.
Hence, the largest gap between any two consecutive $x_i(0)$ must be
at least $1/n$. Thus, there exist some numbers $a$ and $b$, with
$b-a\geq 1/n$, and the set of nodes can be partitioned into two
disjoint subsets $S^-$ and $S_+$ such that $x_i(0)\leq a$ for all
$i\in S_-$, and $x_i(0)\geq b$ for all $i\in S_+$. By Assumption
\ref{as:3}, the graph with arcs $\bigcup_{s=0,\ldots,B-1}
\mathcal{E}(s)$ is connected. Thus, there exists a first time
$\tau\in\{0,1,\ldots,B-1\}$ at which there is a communication
between some node $i\in S_-$ and some node $j\in S_+$, resulting in
an offer from $j$ to $i$. Up until that time,  nodes in $S_-$ have
not interacted with nodes in $S_+$. It follows that $x_k(\tau)\leq
a$ for all $k\in S_-$, and $x_k(\tau)\geq b$ for all $k\in S_+$. In
particular, $x_i(\tau)\leq a$ and $x_j(\tau)\geq b$. There are two
possibilities: either $i$ accepts the offer from $j$, or $i$ accepts
some higher offer, from some other node in $S_+$. In either case, we
conclude that there is a first time $\tau\leq B-1$, at which a node
in $S_-$ accepts an offer from a node in $S_+$.

Let \alexo{us use plain} $x_i$ and $x_j$ \alexo{for} the values at
nodes $i$ and $j$, respectively, at the beginning of period $i$ of
time $\tau$. At the end of that period, the value at both nodes is
equal to $(x_i+x_j)/2$. Thus, the Lyapunov function $V$ decreases by
$$ x_i^2 + x_j^2 - 2 \Big(\frac{x_i + x_j}{2}\Big)^2  =
\frac{1}{2} (x_i - x_j)^2 \geq \frac{1}{2} (b-a)^2\geq  \frac{1}{2n^2}.
$$
At every other time and period, $V$ is nonincreasing, as shown earlier.
Thus, using the inequality $V(0)\leq n$,
$$V(B)\leq V(0)- \frac{1}{2n^2} \leq V(0)\Big(1-\frac{1}{2n^3}\Big).$$
By repeating this argument over the interval $kB,\ldots,(k+1)B$, instead of the interval $0,\ldots,B$, we establish Eq.\ (\ref{eq:vdif}), which concludes the proof.
\end{proof}

\section{Simulations}\label{s:9}

We have proposed several new algorithms for the distributed
consensus and averaging problems. For one of them, namely the
spanning tree heuristic of Section \ref{treeher} (Algorithm
\ref{al:2}), the theoretical performance has been characterized
completely --- see Theorem \ref{th:5} and the discussion at the end
of Section \ref{treeher}.  In this section, we provide simulation
results for the remaining two algorithms.

\subsection{Averaging in fixed networks with two passes of the
agreement algorithm\label{twopassimul}}

\begin{figure}[h]
\begin{center}
$\begin{array}{c@{\hspace{0.1in}}c} \multicolumn{1}{l}{} &
    \multicolumn{1}{l}{} \\ [-0.53cm]
\epsfxsize=2.5in \epsffile{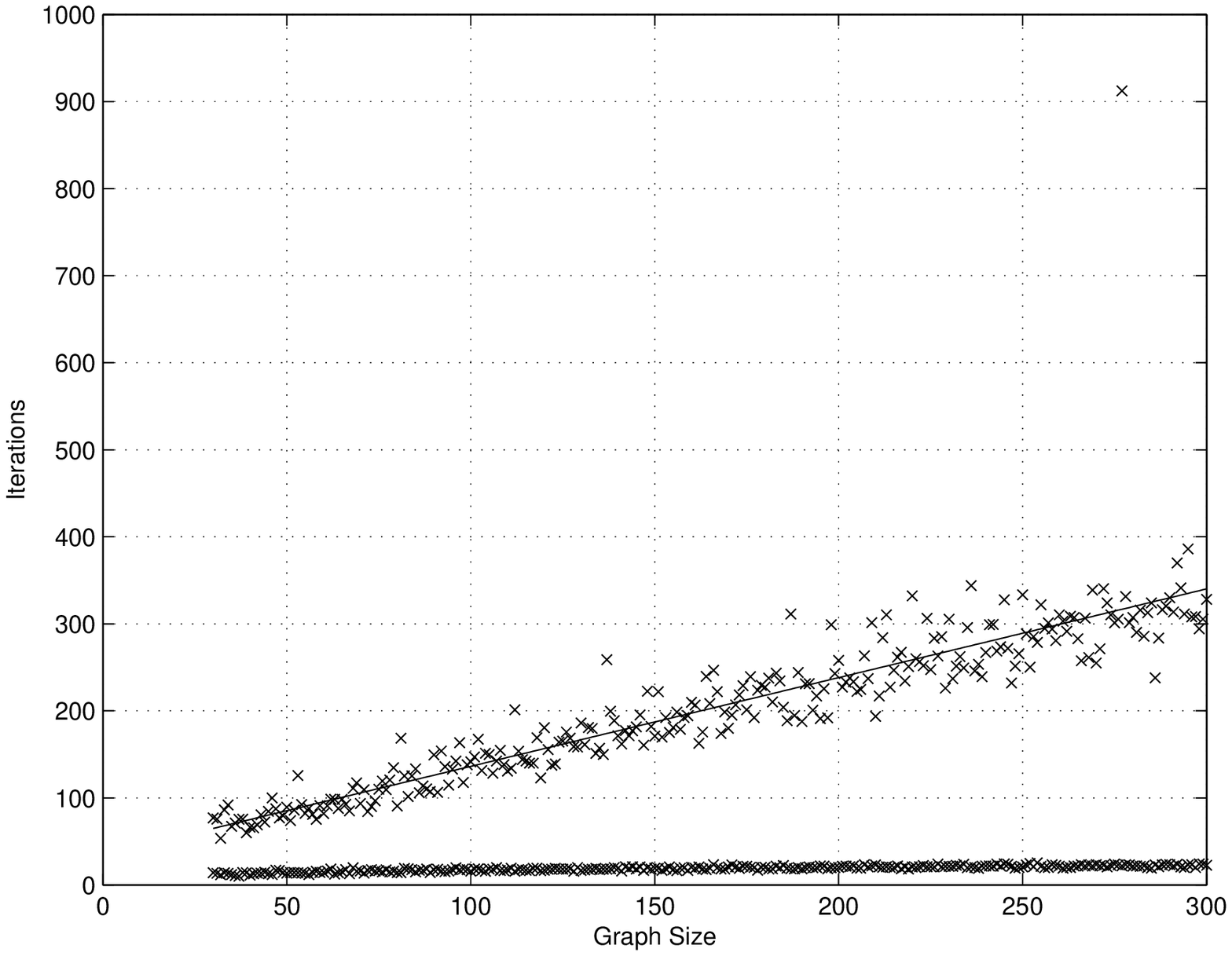} &
    \epsfxsize=2.5in
    \epsffile{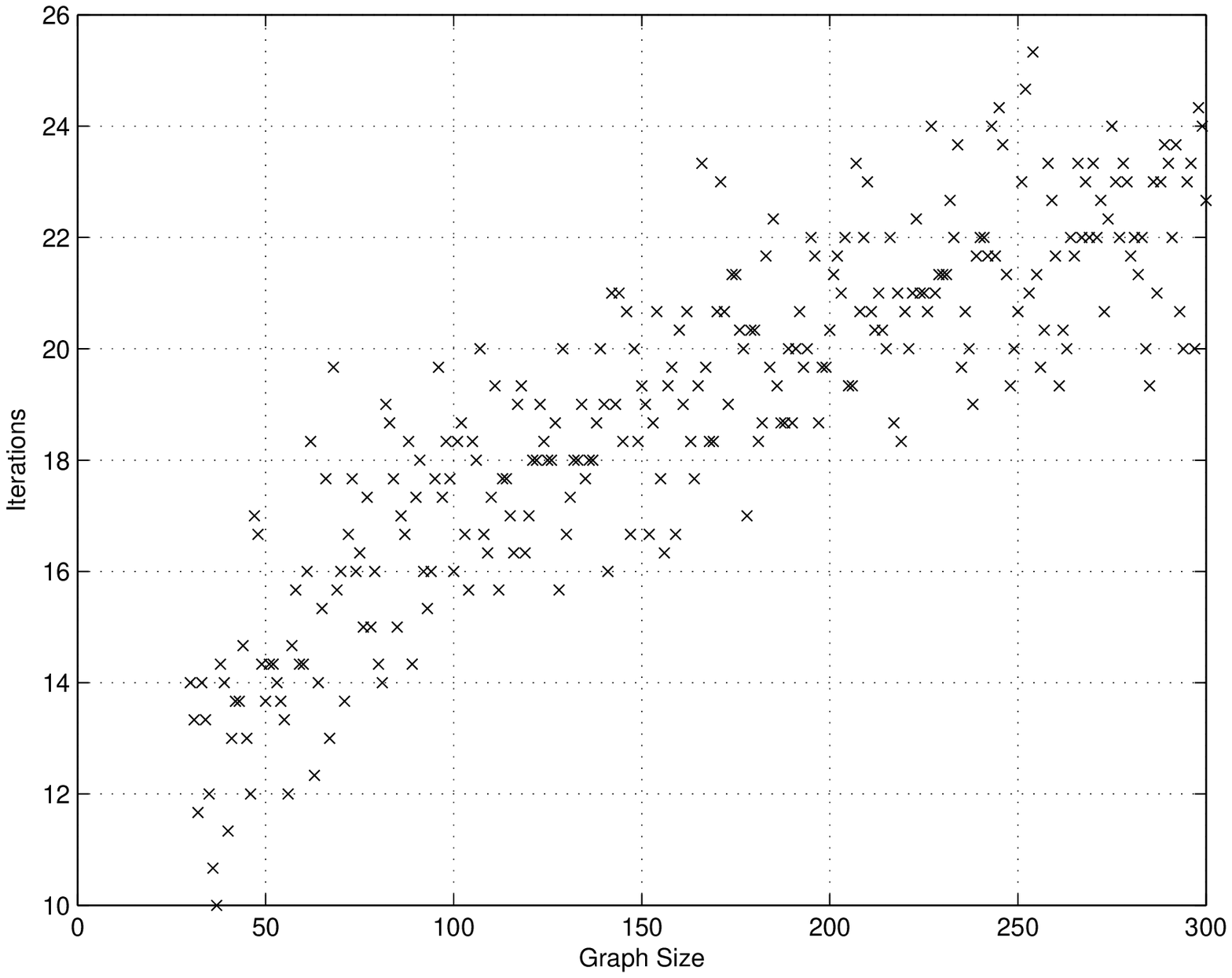} \\ [0.4cm]
\mbox{} & \mbox{}
\end{array}$
\end{center}
\caption{\alexo{On the left:} comparing averaging algorithms on a
geometric random graph. The top line corresponds to the algorithm of
\cite{OSM04}, and the bottom line (close to the horizontal axis)
corresponds to using two parallel passes of the agreement algorithm
(Algorithm \ref{al:1}). \alexo{On the right: a blowup of the
performance of the agreement algorithm.} \label{s1s2}}
\end{figure}

In Section \ref{3a}, we proposed a method for averaging in fixed graphs, based on two parallel executions of the agreement
algorithm (Algorithm \ref{al:1}). We speculated in Section \ref{s:7} that the
presence of a small number of  high degree nodes would make the
performance of our algorithm attractive relative to the algorithm of
\cite{OSM04}, which uses a step size proportional to the inverse of
the largest degree. (Our implementation used a step size of $\epsilon = {1}/{2
d_{\rm max}}$.) Figure \ref{s1s2} presents simulation results for the two
algorithms.

In each simulation, we first generate \alexo{geometric random graph
$G(n,r)$ by placing nodes randomly in $[0,1]^2$ and connecting two
nodes if they are at most $r$ apart. We pick $r=\Theta(\sqrt{\log
n/n})$, which is a standard choice for modeling wireless networks
(c.f. \cite{DSW06}). }

We then change the random graph $G(n,r)$ by picking $n_d$ nodes at
random ($n_d=10$ in both figures) and adding edges randomly to make
the degree of these nodes linear in $n$; this is done by making each
edge incident to these nodes present with probability $1/3$. We run
the algorithm, with random starting values, uniformly distributed in
$[0,1]$, until the largest deviation from the mean is at most
$\epsilon=10^{-3}$.

Each outcome recorded in Figure \ref{s1s2} (for different values of
$n$), is the average of three runs. We conclude that for such
graphs, the convergence time of the algorithm in \cite{OSM04} grows
considerably faster than the one proposed in this paper.

\subsection{Averaging in time-varying Erd\"os-Renyi random graphs}

We report here on simulations involving the load-balancing algorithm
(Algorithm \ref{al:3}) on time-varying random graphs. \alexo{In
contrast to our previous simulations on static geometric graph, we
test two time-varying models which simulate movement. }

\alexo{In both models, we select our initial vector $x(0)$ by
choosing each component independently as a uniform random variable
over $[0,1]$. In our first model,} at each time $t$, we
independently generate an Erd\"os-Renyi random graph $G(t) = G(c,n)$
with $c=3/4$. \alexo{In the second model, at each time step we
independently generate a geometric random graph with $G(n,r)$ with
$r=\sqrt{\log n/n}$. In both models,} if the largest deviation from
the mean is at most $\epsilon = 10^{-3}$, we stop; else, we perform
another iteration of the load-balancing algorithm.

\begin{figure}[h]
\begin{center}
$\begin{array}{c@{\hspace{0.1in}}c} \multicolumn{1}{l}{} &
    \multicolumn{1}{l}{} \\ [-0.53cm]
\epsfxsize=2.5in \epsffile{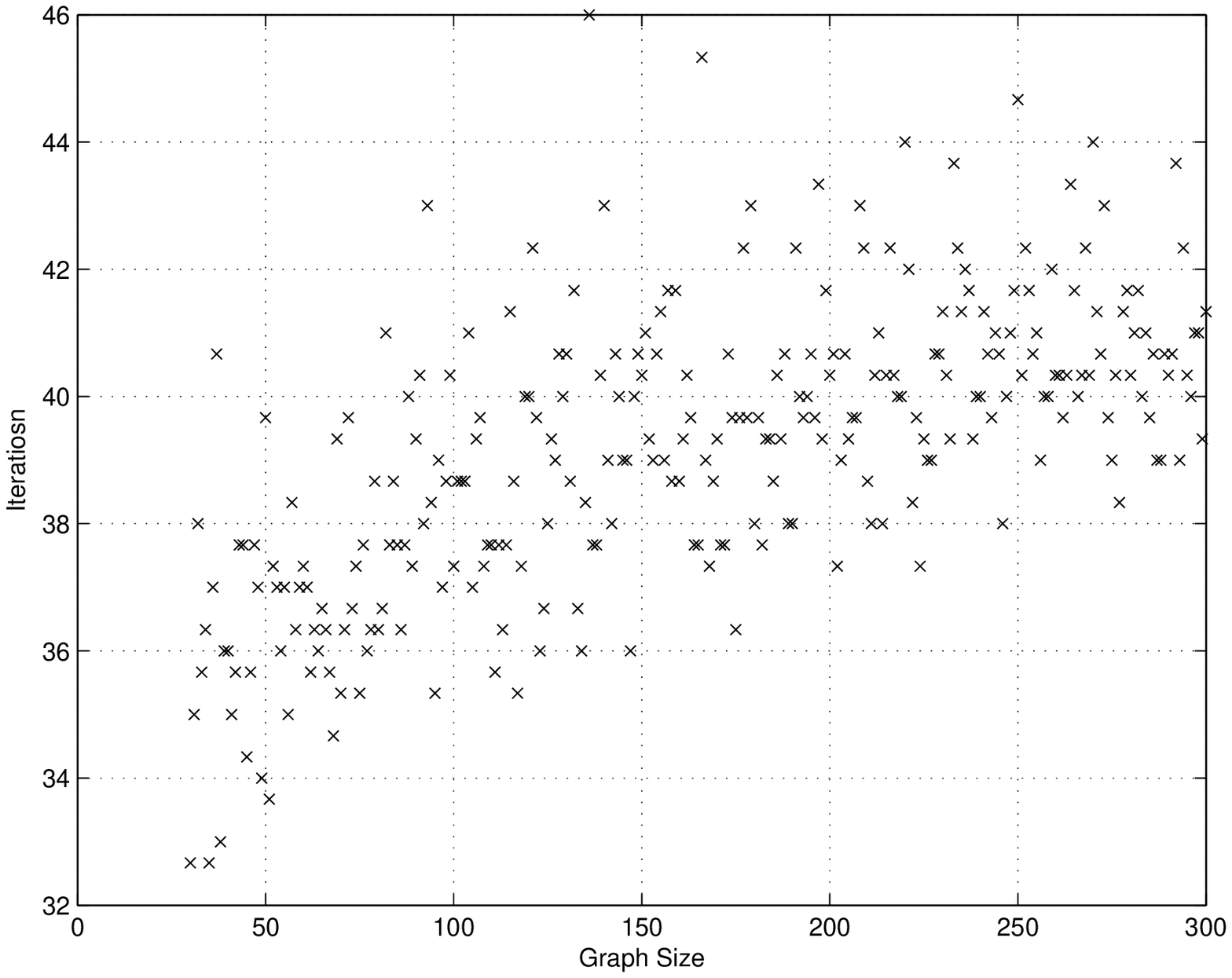} &
    \epsfxsize=2.5in
    \epsffile{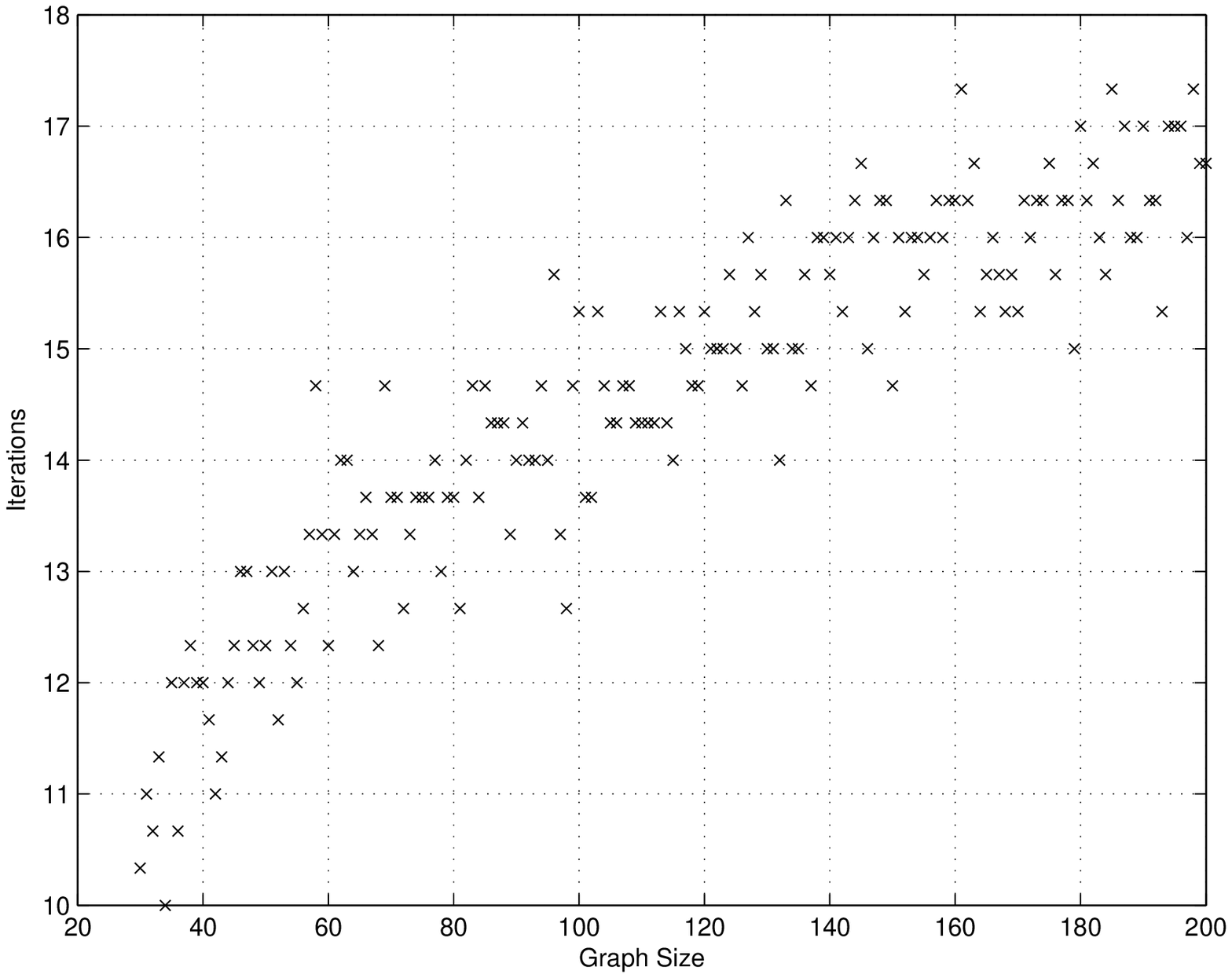} \\ [0.4cm]
\mbox{} & \mbox{}
\end{array}$
\end{center}
\caption{\alexo{On the left:} averaging in time-varying
Erd\"os-Renyi random graphs with the load balancing algorithm. Here
$c=3/4$ at each time $t$. \alexo{On the right:} averaging in
time-varying geometric random graphs with the load balancing
algorithm. Here $r=\sqrt{\log n/n}$.\label{s3s4}}
\end{figure}

The results are summarized in Figure \ref{s3s4}, where again each
point represents the average of three runs. \alexo{We conclude that
in these random models, only a sublinear number of iterations
appears to be needed.}

\section{Concluding remarks}\label{s:10}

In this paper we have considered a variety of consensus and averaging algorithms, and studied their convergence rates. While our discussion was focused on averaging algorithms, several of our results pertain to the closely related consensus problem.

For the case of a fixed topology, we showed that averaging algorithms are easy to construct by using two parallel executions of the agreement algorithm for the consensus problem. We also saw that a reasonable performance guarantee can be obtained by using the equal-neighbor agreement algorithm on a spanning tree, as opposed to a more sophisticated design.

For the case of a fixed topology, the choice of different algorithms is not a purely mathematical issue; one must also take into account the extent to which one is able to design the algorithm offline and provide suitable instructions to each node. After all, if the nodes are able to set up a spanning tree, there are simple distributed algorithms, involving two sweeps along the tree, in opposite directions, with which the sum of their initial values can be computed and disseminated \cite{BT89}, thus eliminating the need for an iterative algorithm.
On the other hand, in less structured environments, with the possibility of occasional changes in the system topology, iterative algorithms can be more resilient. For example, the equal-neighbor agreement algorithm adjusts itself naturally when the topology changes.

In the face of a changing topology (possibly at each time step), the
agreement algorithm continues to work properly, under minimal
assumptions (Theorem \ref{th:1}). On the other hand, its worst-case
convergence time may suffer severely (cf.\ Section \ref{expon}).
Furthermore, it is not apparent how to modify the agreement
algorithm and obtain an averaging algorithm without sacrificing
linearity and/or allowing some additional memory at the nodes. In
Section \ref{s:8}, we introduced an averaging algorithm, which is
nonlinear but leads to a rather favorable (and in particular,
polynomial) convergence time bound. In view of the favorable
performance observed in our simulation results, it would also be
interesting to characterize the average performance of this
algorithm, under a probabilistic mechanism for generating the graphs
$G(t)$, similar to the one in our simulations.

Something to notice about Algorithm \ref{al:3} is that it requires
the topology to remain fixed during each during the exchange of
offers and acceptances/rejections that happens at each step. On the
other hand, without such an assumption, or without introducing a
much larger memory at each node (which would allow for flooding of
individual values), an averaging algorithm may well turn out to be
impossible.

\end{document}